# An Infinite Product for $e^\gamma$ via Hypergeometric Formulas for Euler's Constant $\gamma$

## Jonathan Sondow

**INTRODUCTION.** Euler's constant, $\gamma = 0.5772156649\ldots$, is defined as the limit

$$\gamma = \lim_{n\to\infty}\left(1 + \frac{1}{2} + \cdots + \frac{1}{n} - \ln n\right). \tag{1}$$

In this note we prove the infinite product formula for $e^\gamma = 1.7810724179\ldots$

$$e^\gamma = \left(\frac{2}{1}\right)^{1/2} \left(\frac{2^2}{1\cdot 3}\right)^{1/3} \left(\frac{2^3\cdot 4}{1\cdot 3^3}\right)^{1/4} \left(\frac{2^4\cdot 4^4}{1\cdot 3^6\cdot 5}\right)^{1/5} \cdots, \tag{2}$$

where the $n$th factor is

$$\left(\prod_{k=0}^{n}(k+1)^{(-1)^{k+1}\binom{n}{k}}\right)^{1/(n+1)}.$$

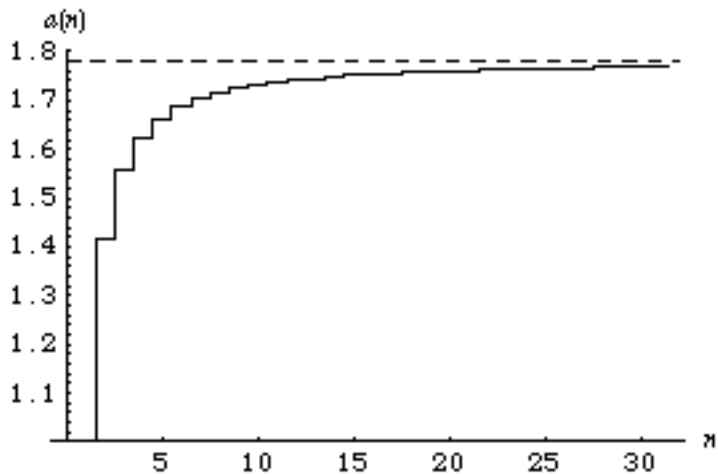

**Figure 1.** The $n$th partial product $a(n)$ of the infinite product for $e^\gamma$

The path that led to (2) began with the proof [5] of the double integral formula for Euler's constant

$$\gamma = \iint_{[0,1]^2} -\frac{1-x}{(1-xy)\ln xy}\,dx\,dy, \tag{3}$$

an analog of double integrals [4] for $\zeta(2)$ and $\zeta(3)$, where $\zeta(s) = \sum_{n\geq 1} n^{-s}$. A series $\sum a_n$ in the proof (essentially the first series in (4) below), with term ratio a quotient of polynomials in $n$, gave the idea to write the series as a hypergeometric function (see [**2**, p. 208]). A hypergeometric transformation then paved the way to (2).

We give four proofs of (2). The first two are hypergeometric; the second uses the digamma function $\psi(t) = \Gamma'(t)/\Gamma(t)$ and is shorter. The last two proofs, which avoid hypergeometrics, depend on integral formulas for $\gamma$: a classical one in the third proof, and (3) in the fourth. Using the classical integral, we end the note with a quick proof of (3).

**PROOF 1.** We follow the latter part of the path described. Using $\ln((n+1)/n) \to 0$ to replace $\ln n$ by $\ln(n+1)$ in (1), we convert the limit to the series

$$\gamma = \sum_{n=1}^{\infty}\left(\frac{1}{n} - \ln\frac{n+1}{n}\right) = \sum_{n=1}^{\infty}\int_1^{\infty}\left(\frac{1}{(t+n-1)^2} - \frac{1}{(t+n-1)(t+n)}\right)dt$$

$$= \sum_{n=1}^{\infty}\int_1^{\infty}\frac{1}{(t+n-1)^2(t+n)}\,dt.$$

Interchanging summation and integration (justified by uniform convergence, the integrand being bounded by $n^{-2}(n+1)^{-1}$), we write the formula as

$$\gamma = \int_1^{\infty}\sum_{n=1}^{\infty}\frac{1}{(t+n-1)^2(t+n)}\,dt = \int_1^{\infty}\frac{1}{t^2(t+1)}\sum_{n=1}^{\infty}\frac{t^2(t+1)}{(t+n-1)^2(t+n)}\,dt$$

$$= \int_1^{\infty}\frac{1}{t^2(t+1)}\left(1 + \frac{t^2}{(t+1)\cdot(t+2)} + \frac{(t(t+1))^2}{(t+1)(t+2)\cdot(t+2)(t+3)} + \cdots\right)dt \tag{4}$$

in order to replace the series by a hypergeometric function [**2**, p. 205], [**3**, p. 101]

$$F\begin{pmatrix}1, & s, & t\\ & u, & v\end{pmatrix} = {}_3F_2\begin{pmatrix}1, & s, & t\\ & u, & v\end{pmatrix}1 = 1 + \frac{s\cdot t}{u\cdot v} + \frac{s(s+1)\cdot t(t+1)}{u(u+1)\cdot v(v+1)} + \cdots, \tag{5}$$

obtaining

$$\gamma = \int_1^{\infty}\frac{1}{t^2(t+1)}\cdot F\begin{pmatrix}1, & t, & t\\ & t+1, & t+2\end{pmatrix}dt.$$

The transformation

$$F\begin{pmatrix}1, & s, & t \\ & s+1, & v\end{pmatrix} = \frac{s}{v-t} \cdot F\begin{pmatrix}1, & v-t, & v-s \\ & v-t+1, & v\end{pmatrix} \tag{6}$$

(valid if $s > 0$ and $v - t > 0$, and proved by applying the identity $F\begin{pmatrix}r, & s, & t \\ & u, & v\end{pmatrix} = F\begin{pmatrix}s, & r, & t \\ & u, & v\end{pmatrix}$, then Thomae's transformation [**3**, pp. 104-105, 111], then the identity) gives the simpler expression

$$\gamma = \int_1^\infty \frac{1}{2t(t+1)} \cdot F\begin{pmatrix}1, & 2, & 2 \\ & 3, & t+2\end{pmatrix} dt. \tag{7}$$

Expanding by (5), and simplifying terms, we get the key formula

$$\gamma = \int_1^\infty \sum_{n=1}^\infty \frac{n!}{(n+1)t(t+1)\cdots(t+n)} dt. \tag{8}$$

Integrating termwise (the $n$th term being bounded by $(n+1)^{-2}$) by the method of partial fractions, we arrive at the series

$$\gamma = \sum_{n=1}^\infty \frac{1}{n+1} \sum_{k=0}^n (-1)^{k+1} \binom{n}{k} \ln(k+1),$$

and exponentiation yields the desired product (2). ●

**PROOF 2.** Using values of the $\psi$ function and a series for $\psi'$ (see [**1**, pp. 15, 16, 22, equations (4), (8), (9), resp.]), we shorten the first proof. From (1) and the formulas $\psi(1) = -\gamma$ and $\psi(n+1) = 1 + 1/2 + \cdots + 1/n - \gamma$, we deduce that

$$\int_1^\infty \left(\psi'(t) - \frac{1}{t}\right) dt = \gamma.$$

Employing (5) and (6), we transform the series $\psi'(t) = \sum_{n\geq 0}(t+n)^{-2}$ into the function

$$\psi'(t) = \frac{1}{t^2} \sum_{n=0}^\infty \frac{t^2}{(t+n)^2} = \frac{1}{t^2} \cdot F\begin{pmatrix}1, & t, & t \\ & t+1, & t+1\end{pmatrix} = \frac{1}{t} \cdot F\begin{pmatrix}1, & 1, & 1 \\ & 2, & t+1\end{pmatrix}.$$

Using (5) (but not (6)), we verify the identity

$$\frac{1}{t}\left[F\begin{pmatrix}1, & 1, & 1 \\ & 2, & t+1\end{pmatrix} - 1\right] = \frac{1}{2t(t+1)} \cdot F\begin{pmatrix}1, & 2, & 2 \\ & 3, & t+2\end{pmatrix},$$

and (7) follows. The rest of the proof is the same. •

*Remark.* One can accelerate the convergence of such hypergeometric representations by computing some terms in (1). For example,

$$\gamma = 1 + \frac{1}{2} + \cdots + \frac{1}{n-1} - \ln n + \int_n^\infty \frac{1}{t}\left[F\left(\begin{matrix}1, & 1, & 1 \\ & 2, & t+1\end{matrix}\right) - 1\right]dt.$$

**PROOF 3.** We reduce the key formula (8) to the classical formula for Euler's constant

$$\gamma = \int_0^1 \left(\frac{1}{\ln u} + \frac{1}{1-u}\right)du. \tag{9}$$

(One way to derive (9) is by inserting the value $\psi(1) = -\gamma$ into Gauss's integral for the $\psi$ function [**6,** p. 246, Example 1].) Denoting integral (8) by $I$, we substitute Euler's beta integral

$$\int_0^1 \frac{u^{t-1}(1-u)^n}{n+1}du = \frac{n!}{(n+1)t(t+1)\cdots(t+n)}$$

(evaluated by $n$-fold integration by parts), then replace $t$ by $t+1$, obtaining

$$I = \int_0^\infty \sum_{n=1}^\infty \int_0^1 \frac{u^t(1-u)^n}{n+1}du\,dt.$$

For fixed $t > 0$, a little calculus shows that

$$\max_{0 \le u \le 1} \frac{u^t(1-u)^n}{n+1} = \frac{1}{n+1} \cdot \frac{n^n t^t}{(n+t)^{n+t}} < \frac{t^t}{n^{1+t}},$$

guaranteeing uniform convergence of the series

$$\sum_{n=1}^\infty \frac{u^t(1-u)^n}{n+1} = \frac{u^t}{1-u}\sum_{n=1}^\infty \frac{(1-u)^{n+1}}{n+1} = u^t\left(\frac{-\ln u}{1-u} - 1\right)$$

on the interval $0 \le u \le 1$. This permits interchanging summation and integration with respect to $u$, which gives

$$I = \int_0^\infty \int_0^1 u^t\left(\frac{-\ln u}{1-u} - 1\right)du\,dt. \tag{10}$$

Reversing the order of integration (the integrand being non-negative), we integrate with respect to $t$, and arrive at integral (9). This proves (8), and (2) follows as before. •

**PROOF 4.** We show that integrals (8) and (3) are equal. Repeating the previous proof through (10), we substitute the integral

$$\frac{1}{1-u}\int_0^{1-u}\frac{v}{1-v}dv = \frac{-\ln u}{1-u}-1,$$

reverse the order of integration from $dv\,du\,dt$ to $dt\,dv\,du$ (non-negativity again), and integrate with respect to *t*, obtaining

$$I = \int_0^1\int_0^{1-u}-\frac{v}{(1-v)(1-u)\ln u}dv\,du. \tag{11}$$

Since the change of variables $v = 1-x, u = xy$ in integral (3) transforms it into integral (11), we infer (8), and (2) follows as usual. •

*Remark.* Integration with respect to *v* in integral (11) produces integral (9) for $\gamma$. Combined with the equality of integrals (3) and (11) by change of variables, this gives a short proof of (3). (The proof in [5] is longer, but uses only the definition of Euler's constant.)

ACKNOWLEDGEMENTS. The author is grateful to P. Sebah and W. Zudilin for suggestions on Thomae's transformation, S. Zlobin for ideas used in the last two proofs, and E. Weisstein for Figure 1.

*209 West 97th Street, New York, NY 10025*
*jsondow@alumni.princeton.edu*